\newcommand{\al}{\alpha}               \newcommand{\be}{\beta}
\newcommand{\ga}{\gamma}               
\newcommand{\de}{\delta}               
\newcommand{\lb}{\lambda}              
\newcommand{\sig}{\sigma}              
               
\newcommand{\veps}{\varepsilon}        \newcommand{\vphi}{\varphi}

\newcommand{\cal}{\mathcal}
           
           \newcommand{\calf}{{\cal F}}

           \newcommand{\calt}{{\cal T}}

\newcommand{\Fix}{{\rm Fix}}

\newcommand{\incl}{\subseteq}        
         
\newcommand{\es}{\emptyset}          
      \newcommand{\limpl}{\Longrightarrow}
  \newcommand{\lequi}{\Longleftrightarrow}
\newcommand{\oo}{\infty}

             \newcommand{\sk}{\smallskip}
       \newcommand{\n}{\noindent}

                \def\R+oo{R_+\cup\{\oo\}}

\def\0dtends   {\stackrel {\it 0d}{\longrightarrow}}
\def\etends   {\stackrel {\it e}{\longrightarrow}}

\newcommand{\barr}{\begin{array}}         \newcommand{\earr}{\end{array}}
\newcommand{\bcor}{\begin{corollary}}     \newcommand{\ecor}{\end{corollary}}
\newcommand{\beq}{\begin{equation}}       \newcommand{\eeq}{\end{equation}}
\newcommand{\bex}{\begin{example}}        \newcommand{\eex}{\end{example}}
\newcommand{\bit}{\begin{itemize}}        \newcommand{\eit}{\end{itemize}}
\newcommand{\blemma}{\begin{lemma}}       \newcommand{\elemma}{\end{lemma}}
\newcommand{\bproof}{\begin{proof}}       \newcommand{\eproof}{\end{proof}}
\newcommand{\bprop}{\begin{proposition}}  \newcommand{\eprop}{\end{proposition}}
\newcommand{\brem}{\begin{remark}}        \newcommand{\erem}{\end{remark}}
\newcommand{\btab}{\begin{tabular}}       \newcommand{\etab}{\end{tabular}}
\newcommand{\btheorem}{\begin{theorem}}   \newcommand{\etheorem}{\end{theorem}}

\documentclass[reqno]{amsart}

\newtheorem{theorem}{\bf Theorem}
\newtheorem{corollary}{\bf Corollary}
\newtheorem{example}{\bf Example}
\newtheorem{lemma}{\bf Lemma}
\newtheorem{proposition}{\bf Proposition}
\newtheorem{remark}{\bf Remark}

\begin{document}

\title
[Function contractive maps in partial metric spaces]
{FUNCTION CONTRACTIVE MAPS IN \\
PARTIAL METRIC SPACES}

\author{Mihai Turinici}
\address{
"A. Myller" Mathematical Seminar;
"A. I. Cuza" University;
700506 Ia\c{s}i, Romania
}
\email{mturi@uaic.ro}




\begin{abstract}
Some fixed point results are given
for a class of functional contractions
over partial metric spaces. 
These extend some contributions in the area due to
Ili\'{c} et al [Math. Comput. Modelling, 55 (2012), 801-809]. 
\end{abstract}

\maketitle

\section{Introduction}
\setcounter{equation}{0}

Let $(X,d)$ be a metric space; 
and $T:X\to X$ be a selfmap of $X$.
Denote the class of  its fixed points in $X$
as $\Fix(T)=\{z\in X; z=Tz\}$.
Concerning the existence and uniqueness of such 
elements, a basic result is the 
1922 one due to 
Banach \cite{banach-1922};
it says that, if $(X,d)$ is complete and 
\bit
\item[(a01)]
$d(Tx,Ty)\le \lb d(x,y)$,\ \ $\forall x,y\in X$,
\eit
for some $\lb$ in $[0,1[$ then, the following conclusions hold:

{\bf i)} $\Fix(T)$ is a singleton, $\{z\}$;\ \ 
{\bf ii)} $T^nx\to z$ as $n\to \oo$, for each $x\in X$.

\n
This result found a multitude of applications in 
operator equations theory;
so, it was the subject of many extensions.
For example, a natural way of doing this  
is by considering "functional" contractive conditions
of the form
\bit
\item[(a02)]
$d(Tx,Ty)\le F(d(x,y),d(x,Tx),d(y,Ty),d(x,Ty),d(y,Tx))$,
\ \ $\forall x,y\in X$;
\eit
where $F:R_+^5\to R_+$ is an appropriate function.
For more details about the possible choices of $F$
we refer to 
Collaco and E Silva \cite{collaco-e-silva-1997},
Kincses and Totik \cite{kincses-totik-1990}.
Park \cite{park-1980}
and
Rhoades \cite{rhoades-1977};
see also 
Turinici \cite{turinici-1976}.
Another way of extension is that of conditions imposed
upon $d$ being modified.
Some results of this type were obtained in the last 
decade over the realm of partial metric spaces,
introduced as in 
Matthews \cite{matthews-1994}.
However, 
for the subclass of zero-complete partial metric spaces, 
most of these statements 
based on contractive conditions like in (a02)
are nothing but clones of their corresponding 
(standard) metrical counterparts;
see, for instance,
Abdeljawad  \cite{abdeljawad-2011},
Altun et al \cite{altun-sola-simsek-2010},
Altun and Sadarangani \cite{altun-sadarangani-2010},
Duki\'{c} et al \cite{dukic-kadelburg-radenovic-2011},
Karapinar and Erhan \cite{karapinar-erhan-2011},
Oltra and Valero  \cite{oltra-valero-2004},
Paesano and Vetro \cite{paesano-vetro-2012},
Romaguera \cite{romaguera-2012},
Valero \cite{valero-2005},
and the references therein.
Concerning the non-reducible to (a02) contractive 
maps in complete partial metric spaces, 
the first specific result in the area was obtained by
Ili\'{c} et al \cite{ilic-pavlovic-rakocevic-2011}.
It is our aim in the following to show that 
all these are again reducible to a line of argument
used in the standard metrical case.
Further aspects will be delineated elsewhere.

\section{Partial metrics}
\setcounter{equation}{0}

Let $X$ be a nonempty set.
By a {\it symmetric} on $X$ we mean any map $d:X\times X\to R_+$
with
\bit
\item[(b01)]
$d(x,y)=d(y,x)$,\ \ $\forall x,y\in X$.
\eit
Fix such an object; and put (for $x,y\in X$)
\bit
\item[]
$b(x,y)=(1/2)(d(x,x)+d(y,y))$, \ \
$c(x,y)=\max\{d(x,x),d(y,y)\}$. 
\eit
Call the (symmetric) $d$, {\it reflexive-triangular}, if
\bit
\item[(b02)]
$d(x,z)\le d(x,y)+d(y,z)-d(y,y)$,\ \ $\forall x,y,z\in X$.
\eit
According to
Matthews \cite{matthews-1994},
we say that the symmetric $d$ is a {\it partial metric}
when it is reflexive-triangular and
\bit
\item[(b03)]
$(b(x,y)\le) c(x,y)\le d(x,y)$,\ \ $\forall x,y\in X$
\hfill (Matthews property)
\item[(b04)]
$d(x,y)=d(x,x)=d(y,y)$ $\limpl$ $x=y$
\hfill (weak sufficiency).
\eit
In this case, $(X,d)$ is called a {\it partial metric space}.
Denote in what follows
\bit
\item[(b05)]
$e(x,y)=2(d(x,y)-b(x,y))$,\ \ $x,y\in X$.
\eit

\blemma \label{le1}
The mapping $e(.,.)$ is a (standard) metric on $X$;
in addition,
\beq \label{201}
|d(x,x)-d(y,y)|\le e(x,y),\ \ \forall x,y\in X.
\eeq
\elemma

\bproof
The first part is clear; so, we do not give details.
For the second one, let $x,y\in X$ be arbitrary fixed.
Without loss, one may assume that $d(x,x)\ge d(y,y)$; 
for, otherwise, we simply interchange $x$ and $y$.
In this case, (\ref{201}) becomes 
$$
d(x,x)-d(y,y)+d(x,x)+d(y,y)\le 2d(x,y).
$$
This is evident, from the Matthews property; hence the conclusion.
\eproof

Technically speaking, all constructions and results
below are ultimately based on the (standard) 
metric $e$; but, for an appropriate handling, 
these are expressed in terms of the partial metric $d$.
\sk

{\bf (A)}
Denote, for $x\in X$, $\veps> 0$,
\bit
\item[]
$X_d(x,\veps)=\{y\in X; d(x,y)< d(x,x)+\veps\}$
\eit
(the open $d$-sphere with center $x$ and radius $\veps$).
Let $\calt(d)$ stand for the topology having as basis the family 
$\{X_d(x,\veps); x\in X, \veps> 0\}$; it is $T_0$-separated,
as established in
Bukatin et al \cite{bukatin-kopperman-matthews-pajoohesh-2009}. 
The sequential convergence structure attached to this topology
may be depicted as follows.
Let $N:=\{0,1,...\}$ stand for the set of all natural numbers;
and, for each $k\in N$, put $N(k,\le)=\{n\in N; k\le n\}$.
By a sequence over (the nonempty set) $A$, 
we shall mean any map $n\mapsto x(n):=x_n$ from $N(n_0,\le)$ to $A$,
where $n_0$ depends on this map only;
it will be also written as $(x_n; n\ge n_0)$, or $(x_n)$,
when $n_0$ is generically understood. 
Likewise, by a doubly indexed sequence in $A$ we shall mean any map
$(m,n)\mapsto y(m,n):=y_{m,n}$ from 
$N(p_0,\le)\times N(p_0,\le)$ 
to $A$, where $p_0$ depends on this map only;
it will be also written as $(y_{m,n}; m,n\ge p_0)$,
or $(y_{m,n})$, where $p_0$ is generically understood.
Now, for the sequence $(x_n; n\ge n_0)$ in $X$
and the point $x\in X$, define
\bit
\item[]
$(x_n; n\ge n_0)\to x$ (modulo $\calt(d)$) iff \\
$\forall \veps> 0$, $\exists n(\veps)\ge n_0$:\ 
$n\ge n(\veps)$ $\limpl$ $x_n\in X_d(x,\veps)$.
\eit
When $n_0$ is generic, it will be also written 
as $x_n\to x$ (modulo $\calt(d)$).
Note that, by the very definition of our topology,
\beq \label{202}
\mbox{
$x_n\to x$ (modulo $\calt(d)$) iff
$\lim_n d(x_n,x)=d(x,x)$.
}
\eeq
Further, let $\calt(e)$ stand for the 
Hausdorff topology attached to the metric $e$.
We shall write $(x_n; n\ge n_0)\to x$ (modulo $\calt(e)$)
as $(x_n; n\ge n_0) \etends x$; or, simply, $x_n\etends x$. 
By definition, this means 
$e(x_n,x)\to 0$ as $n\to \oo$;
i.e.: 
$$
\forall \al> 0, \exists n(\al)\ge n_0:\ \ 
n\ge n(\al)\ \limpl e(x_n,x)< \al.
$$
A useful interpretation of this in terms of $d$
may be given as follows.

\blemma \label{le2}
Under the above conventions,
\beq \label{203}
\mbox{
$x_n\etends x$ iff 
$\lim_n d(x_n,x)=\lim_{m,n} d(x_m,x_n)=d(x,x)$.
}
\eeq
\elemma

\bproof
If (\ref{203}) holds, then
$$
\lim_nd(x_n,x)=d(x,x),\ \ \lim_n d(x_n,x_n)=d(x,x);
$$
so that,
$$
\lim_n e(x_n,x)=\lim_n[2d(x_n,x)-d(x_n,x_n)-d(x,x)]=0.
$$
Conversely, assume that $e(x_n,x)\to 0$ as $n\to \oo$.
By definition, this means
$$
[d(x_n,x)-d(x_n,x_n]+[d(x_n,x)-d(x,x)]\to 0.
$$
Since the quantities in the square brackets are positive, 
this yields
$\lim_nd(x_n,x)=\lim_nd(x_n,x_n)=d(x,x)$.
Finally, let $m,n$ be arbitrary fixed.
By the reflexive-triangular property,
$$
d(x_m,x_n)\le d(x_m,x)+d(x,x_n)-d(x,x),
$$
$$
d(x_m,x_n)\ge d(x_m,x)-d(x,x_n)+d(x_n,x_n).
$$
Since the limit in the right hand side 
is $d(x,x)$ when $m,n\to \oo$, we are done.
\eproof

Note that, as a consequence of this, the 
convergence structures attached to $\calt(d)$
and $\calt(e)$ are quite distinct; 
we do not give details.
\sk

{\bf (B)}
Remember that, a sequence $(x_n; n\ge n_0)$ in $X$ is 
$e$-Cauchy when $e(x_m,x_n)\to 0$ as $m,n\to \oo$;
that is,
$$
\forall \al> 0, \exists p(\al)\ge n_0:\ \ 
m,n\ge p(\al)\ \limpl e(x_m,x_n)< \al.
$$
As before, we are interested to 
characterize this property in terms of $d$.

\blemma \label{le3}
The generic property below is true: 
\beq \label{204}
\mbox{
$(x_n; n\ge n_0)$\ is 
$e$-Cauchy $\lequi$ 
$\lim_{m,n}d(x_m,x_n)$ exists (in $R_+$).
}
\eeq
\elemma

\bproof
i) 
Assume that $(x_n; n\ge n_0)$ 
is taken as in the right hand of (\ref{204}); i.e.,
$$
\mbox{
$\lim_{m,n} d(x_m,x_n)=\ga$,\ \ for some $\ga\in R_+$.
}
$$
In particular, (for $m=n$), 
we must have $\lim_nd(x_n,x_n)=\ga$. But then, 
$$
\lim_{m,n}e(x_m,x_n)=
\lim_{m,n}[2d(x_m,x_n)-d(x_m,x_m)-d(x_n,x_n)]=0;
$$
so that, $(x_n)$ is $e$-Cauchy.

ii)
Assume that $(x_n; n\ge n_0)$ is $e$-Cauchy; 
that is: $\lim_{m,n}e(x_m,x_n)=0$.
By the very definition of our metric, this means
\beq \label{205}
\lim_{m,n}[d(x_m,x_n)-d(x_m,x_m)+d(x_m,x_n)-d(x_n,x_n)]=0.
\eeq
Let $\al> 0$ be arbitrary fixed.
By the imposed hypotheses and Lemma \ref{le1},
there exists $k:=k(\veps)\ge n_0$ such that
$$
|d(x_m,x_m)-d(x_n,x_n)|\le e(x_m,x_n)< \al,\ \ 
\forall m,n\ge k.
$$
This tells us that $(d(x_n,x_n); n\ge n_0)$ is a 
Cauchy sequence in $R_+$; wherefrom 
$\lim_n d(x_n,x_n)=\ga$, for some $\ga\in R_+$.
Combining with (\ref{205}),
$\lim_{m,n} d(x_m,x_n)=\ga$; 
hence, $(x_n)$ fulfills the property 
of the right hand side in (\ref{204}).
\eproof

{\bf (C)}
Denote, for each sequence $(x_n; n\ge n_0)$ in $X$,
\bit
\item[]
$e-\lim_n(x_n)=\{x\in X; x_n\etends x\}$.
\eit
As $e(.,.)$ is a metric, $e-\lim_n(x_n)$ is either 
empty or a singleton;
in this last case, 
we say that $(x_n)$ is {\it $e$-convergent}.
Again by the metric property of $e$, we have that
each $e$-convergent sequence 
is $e$-Cauchy too; when the reciprocal holds, 
$(X,e)$ is referred to as {\it complete}.

In the following, a useful result is given
about the $e$-semi-Cauchy sequences that are not
endowed with the $e$-Cauchy property.
Some conventions are needed.
Given a sequence $(\tau_n; n\ge n_0)$ in $R$ 
and a point $\tau\in R$, we have the equivalence 
$$
\mbox{
$\tau_n\to \tau$ iff $(\tau_n; n\ge k) \to \tau$, for all $k\ge n_0$.
}
$$
Also, for the same initial data, put
\bit
\item[]
$\tau_n \downarrow \tau$ if 
[$\tau_n\ge \tau, \forall n\ge n_0$] and $\tau_n\to \tau$
\item[]
$\tau_n \searrow \tau$ if 
$(\tau_n; n\ge k) \downarrow \tau$, for some $k\ge n_0$.
\eit
The implication below is evident:
$$
\tau_n \downarrow \tau \limpl 
(\tau_n; n\ge k) \downarrow \tau,\ \ \forall k\ge n_0.
$$
On the other hand, if $(\tau_n; n\ge n_0) \searrow \tau$ 
then $\tau_n \downarrow \tau$ is not in general true.
Likewise, given the double indexed sequence
$(\sig_{m,n}; m,n\ge p_0)$ in $R$ and the point 
$\sig\in R$, 
$$
\mbox{
$\sig_{m,n}\to \sig$ iff $(\sig_{m,n}; m,n\ge k) \to \sig$, for all $k\ge p_0$.
}
$$
Also, for the same initial data, put
\bit
\item[]
$\sig_{m,n} \downarrow \sig$ if 
[$\sig_{m,n}\ge \sig, \forall m,n\ge p_0$] and $\sig_{m,n}\to \sig$
\item[]
$\sig_{m,n} \searrow \sig$ if 
$(\sig_{m,n}; m,n\ge k) \downarrow \sig$, for some $k\ge 0$.
\eit
As before, the implication below is evident:
$$
\sig_{m,n} \downarrow \sig \limpl 
(\sig_{m,n}; m,n\ge k) \downarrow \sig,\ \ \forall k\ge p_0.
$$
On the other hand, if $\sig_{m,n} \searrow \sig$ 
then $\sig_{m,n} \downarrow \sig$ is not in general true.

Now, given a sequence $(x_n; n\ge n_0)$ in $X$, 
call it {\it $e$-semi-Cauchy} provided
\bit
\item[]
$(\al_n:=d(x_n,x_n)) \searrow \ga$,
$(\rho_n:=d(x_n,x_{n+1})) \searrow \ga$,\ \ for some $\ga\in R_+$.
\eit
Note that, in such a case, the $e$-Cauchy property of
$(x_n)$ means $d(x_m,x_n)\searrow \ga$.

\bprop \label{p1}
Suppose that $(x_n; n\ge n_0)$ is $e$-semi-Cauchy but not $e$-Cauchy.
There exists then $k\ge n_0$, $\veps> 0$, $j(\veps)\ge k$, and a couple
of rank-sequences $(m(j); j\ge k)$, $(n(j); j\ge k)$ with
\beq \label{206}
j\le m(j)\le n(j),\ d(x_{m(j)},x_{n(j)})\ge \ga+\veps,\ 
\forall j\ge k
\eeq
\beq \label{207}
n(j)-m(j)\ge 2,\ d(x_{m(j)},x_{n(j)-1})< \ga+\veps,\ 
\forall j\ge j(\veps)
\eeq
\beq \label{208}
(d(x_{m(j)},x_{n(j)}); j\ge k) \downarrow \ga+\veps
\eeq
\beq \label{209}
(d(x_{m(j)+p},x_{n(j)+q}); j\ge k) \to \ga+\veps,\ 
\forall p,q\in \{0,1\}.
\eeq
\eprop

\bproof
By the very definition of $e$-semi-Cauchy property,
there must be some $k\ge n_0$ in such a way that
$(\al_n; n\ge k)\downarrow \ga$, 
$(\rho_n; n\ge k)\downarrow \ga$.
In this case, the conclusion to be derived writes
$(d(x_m,x_n); m,n\ge k) \downarrow \ga$; i.e.,
$$
\forall \veps> 0,\ \exists j=j(\veps)\ge k:\ 
j\le m\le n \limpl d(x_m,x_n)< \ga+\veps.
$$
The negation of this property means
that, there exists $\veps> 0$, such that:
for each $j\ge k$, there may be found a couple of
ranks $m(j)$, $n(j)$ with the property (\ref{206}).
Given this $\veps> 0$, there exists, 
by the $e$-semi-Cauchy property, some $j(\veps)\ge k$ with
\beq \label{210}
\ga \le \al_i\le \rho_i< \ga+\veps,\ \ \forall i\ge j(\veps).
\eeq
We claim that this is our desired rank for the remaining
conclusions in the statement.
In fact, for each $j\ge j(\veps)$, let $n(j)$ 
be the minimal rank fulfilling (\ref{206})
(with $n$ in place of $n(j)$).
For the moment, (\ref{207}) is clear, via (\ref{210}).
Further, by this relation and the reflexive-triangular inequality,
$$  \barr{l}
\ga+\veps\le d(x_{m(j)},x_{n(j)})\le 
d(x_{m(j)},x_{n(j)-1})+\rho_{n(j)-1}-\al_{n(j)-1} \le \\
\ga+\veps+\rho_{n(j)-1}-\al_{n(j)-1},\ \forall j\ge j(\veps).
\earr
$$
So, passing to limit as $j\to \oo$ yields (\ref{208}).
Finally, again by the reflexive-triangular inequality,
one has, for all $j\ge j(\veps)$,
$$
d(x_{m(j)},x_{n(j)+1})\le 
d(x_{m(j)},x_{n(j)})+\rho_{n(j)}-\al_{n(j)},
$$
$$
d(x_{m(j)},x_{n(j)+1})\ge 
d(x_{m(j)},x_{n(j)})-\rho_{n(j)}+\al_{n(j)+1}.
$$
By a limit process upon $j$  one gets 
the case $(p=0,q=1)$ of (\ref{209}).
The remaining ones may be obtained in a similar way.
\eproof

\section{Normal functions}
\setcounter{equation}{0}

Let $\calf(A)$ stand for the class of all functions from 
$A$ to itself.

{\bf (A)}
Call $\vphi\in \calf(R_+)$, {\it normal} when
\bit
\item[(c01)]
$\vphi(0)=0$; $\vphi(t)< t$,\ for all $t$ in $R_+^0:=]0,\oo[$.
\eit
For any such function $\vphi$, and any $s$ in $R_+^0$,
denote
\bit
\item[(c02)]
$\limsup_{t\to s+}\vphi(t)=\inf_{\veps> 0}\Phi[s+](\veps)$, where \\
$\Phi[s+](\veps)=\sup\{\vphi(t); s\le t< s+\veps\}$, $\veps> 0$ 
\item[(c03)]
$\limsup_{t\to s}\vphi(t)=\inf_{\veps> 0}\Phi[s](\veps)$, where \\
$\Phi[s](\veps)=\sup\{\vphi(t); s-\veps< t< s+\veps\}$, $\veps> 0$. 
\eit
Note that, again by the choice of $\vphi$, one has
\beq \label{301}
\vphi(s)\le\limsup_{t\to s_+}\vphi(t)\le 
\limsup_{t\to s}\vphi(t)\le s,\ \ \forall s\in R_+^0.
\eeq

\blemma \label{le4}
Let $\vphi\in \calf(R_+)$ be normal;
and $s\in R_+^0$ be arbitrary fixed. Then,
 
{\bf i)} 
$\limsup_n \vphi(t_n)\le \limsup_{t\to s_+}\vphi(t)$,
for each sequence $(t_n; n\ge n_0)$ with $t_n\searrow s$

{\bf ii)} 
$\limsup_n \vphi(t_n)\le \limsup_{t\to s}\vphi(t)$,
for each sequence $(t_n; n\ge n_0)$ with $t_n\to s$.
\elemma

\bproof
By definition, there exists $k\ge n_0$ with 
$(t_n; n\ge k)\downarrow s$.
Given $\veps> 0$, 
there exists a rank $p(\veps)\ge k$ such that
$s\le t_n< s+\veps$, $\forall n\ge p(\veps)$; whence
$$
\limsup_n \vphi(t_n)\le \sup\{\vphi(t_n); n\ge p(\veps)\}\le 
\Phi[s+](\veps).
$$
It suffices taking the infimum (=limit) as $\veps\to 0+$ 
in this relation to get the desired fact.
The second part is proved in a similar way.
\eproof

We say that the normal function $\vphi\in \calf(R_+)$ is 
{\it right limit normal} (respectively: {\it limit normal}), if 
\bit
\item[(c04)]
$\limsup_{t\to s_+}\vphi(t)< s$\ 
(respectively: $\limsup_{t\to s}\vphi(t)< s$), \ \ 
$\forall s\in R_+^0$.
\eit
In particular, for the normal function $\vphi\in \calf(R_+)$, 
this holds whenever $\vphi$ is right usc (respectively: usc) on $R_+^0$. 
Note that this property is fulfilled when 
$\vphi$ is right continuous 
(respectively: continuous) on $R_+^0$.
\sk

{\bf (B)}
Let $\psi\in \calf(R_+)$ be a function. Denote
\bit
\item[(c05)]
$\liminf_{t\to \oo}\psi(t)=\sup_{\al\ge 0}\Psi(\al)$,\
where $\Psi(\al)=\inf\{\psi(t); t\ge \al\}$, $\al\ge 0$.
\eit
Call $\psi$, {\it semi-coercive}, provided 
$\lim\inf_{t\to \oo} \psi(t)> 0$.
For example, this is valid whenever $\psi$ is {\it coercive}:
$\lim_{t\to \oo} \psi(t)=\oo$.

\blemma \label{le5}
Let $\psi\in \calf(R_+)$ be some function; 
and $(t_n; n\ge n_0)$ be a sequence with
$(t_n; n\ge n_0)\to \oo$. Then,

{\bf iii)}
$\liminf_n\psi(t_n)\ge \liminf_{t\to \oo}\psi(t)$, 

{\bf iv)}
$\liminf_n\psi(t_n)> 0$, whenever $\psi$ is semi-coercive.
\elemma

The proof is very similar to the one of Lemma \ref{le4};
so, we omit it.

{\bf (C)}
Let $(a_n; n\ge n_0)$ be a bounded in $R_+$ sequence.
From this choice, 
$A_n:=\sup\{a_h; h\ge n\}$ exists in $R_+$, for all $n\ge n_0$;
moreover, $(A_n; n\ge n_0)$ is descending.
Now, by definition, 
$\limsup_n(a_n)=\inf_{n\ge n_0}(A_n)$; 
this, by the descending property of this last sequence, writes:
$\limsup_n(a_n)=\lim_n(A_n)$.

\blemma \label{le6}
Let $F:R_+^3\to R_+$ be a function with
\bit
\item[(c06)]
$F$ is increasing in all variables,
\item[(c07)]
$F$ is continuous at the right over $R_+^3$;
\eit
and let
$(a_n; n\ge n_0)$, $(b_n; n\ge n_0)$, $(c_n; n\ge n_0)$
 be bounded sequences. Then, 
$$
\limsup_n F(a_n,b_n,c_n)\le F(a,b,c),
$$
where 
$a=\limsup_n(a_n)$, $b=\limsup_n(b_n)$, $c=\limsup_n(c_n)$.
\elemma

\bproof
Denote (for each $n\ge n_0$)
$A_n:=\sup\{a_h; h\ge n\}$, 
$B_n:=\sup\{b_h; h\ge n\}$,
$C_n:=\sup\{c_h; h\ge n\}$.
Let $\veps> 0$ be arbitrary fixed.
By the remark above, there exists $h(\veps)\ge n_0$
in such a way that
$A_n< a+\veps$, $B_n< b+\veps$, $C_n< c+\veps$, 
for all $n\ge h(\veps)$.
This (by the properties of $F$), yields for all $n\ge h(\veps)$,
$$
\limsup_n F(a_n,b_n,c_n)\le 
\sup\{F(a_n,b_n,c_n); n\ge h(\veps)\}\le F(a+\veps,b+\veps,c+\veps).
$$
Passing to limit as $\veps\to 0$ gives the desired fact.
\eproof

As we shall see, the usual particular case to be considered is that of
$F(a,b,c)=a+\max\{b,c\}$, $a,b,c\in R_+$.
Clearly, (c06) and (c07) are fulfilled here.

\section{Main result}
\setcounter{equation}{0}

Let $X$ be a nonempty set.
Take a partial metric $d(.,.)$ over it; 
and let $e(.,.)$ stand for its associated
(standard) metric. Assume in the following that
\bit
\item[(d01)]
$(X,e)$ is complete.
\eit
Let $T:X\to X$ be a selfmap of $X$.
We say that $z\in X$ is a {\it $d$-fixed point}
of $T$, when $d(z,Tz)=d(z,z)$.
The set of all these points will be denoted as
$\Fix(T;d)$.

Loosely speaking, such points are to be determined 
as $e$--limits of iterative processes.
Precisely, let us say that $x\in X$ is a
{\it $(d;T)$-Picard point}, when
\bit
\item[(d02)]
$(T^nx; n\ge 0)$ is $e$-convergent and 
$x^*:=e-\lim_n(T^nx)$ is in $\Fix(T;d)$.
\eit
If this holds for each $x\in X$, then $T$ 
will be referred to as a {\it $d$-Picard operator};
cf. Rus \cite[Ch 2, Sect 2.2]{rus-2001}.

Now, sufficient conditions for getting this property
are of functional contractive type.
Denote, for $x,y\in X$,
\bit
\item[]
$M_1(x,y)=\max\{d(x,y),d(x,Tx),d(y,Ty)\}$,\
$M_2(x,y)=\\
(1/2)[d(x,Ty)+d(Tx,y)]$,\ 
$M_3(x,y)=\max\{M_1(x,y),M_2(x,y)\}$.
\eit
Note that, if $x=y$, then (by the Matthews property)
$$
M_1(x,x)=\max\{d(x,x),d(x,Tx)\}=d(x,Tx),\
M_2(x,x)=d(x,Tx);
$$
hence
\beq \label{401}
M(x,x)=d(x,Tx),\ \ x\in X.
\eeq
Likewise, if $y=Tx$, then
$$ \barr{l}
M_1(x,Tx)=\max\{d(x,Tx),d(Tx,T^2x)\}, \\
M_2(x,Tx)=(1/2)[d(x,T^2x)+d(Tx,Tx)]\le \\
(1/2)[d(x,Tx)+d(Tx,T^2x)]\le \max\{d(x,Tx),d(Tx,T^2x)\};
\earr
$$
hence 
\beq \label{402}
M(x,Tx)=\max\{d(x,Tx),d(Tx,T^2x)\},\ \  x\in X.
\eeq
Given the normal function $\vphi\in \calf(R_+)$ 
and $g\in \{b,c\}$, let us say that $T$ is 
{\it $(M;g;\vphi)$-contractive} (modulo $d$) when
\bit
\item[(d03)]
$d(Tx,Ty)\le 
\max\{\vphi(M(x,y)),g(x,y)\},\ \ \forall x,y\in X$.
\eit
We are now in position to give the first main result of this
exposition.

\btheorem \label{t1}
Suppose that $T$ is $(M,c;\vphi)$-contractive (modulo $d$)
for some right limit normal function $\vphi\in \calf(R_+)$.
Then, each $x\in X$ is a $(d;T)$-Picard point, in the sense

{\bf i)}
$(\rho_n(x):=d(T^nx,T^{n+1}x); n\ge 0)$ is descending;
hence, $\rho_n(x)\downarrow \ga(x)$, for some $\ga(x)\in R_+$

{\bf ii)}
$(\al_n(x):=d(T^nx,T^nx); n\ge 0) \searrow \ga(x)$;
so that, $(T^nx; n\ge 0)$ is $e$-semi-Cauchy

{\bf iii)}
$(\de_{m,n}(x):=d(T^mx,T^nx); m,n\ge 0) \searrow \ga(x)$;
whence, $(T^nx; n\ge 0)$ is $e$-Cauchy

{\bf iv)}
$x^*:=e-\lim_n(T^nx)$ is an element of $\Fix(T;d)$.

\n
Hence, in particular, $T$ is a $d$-Picard operator.
\etheorem

\bproof
Fix in the following some point $x_0\in X$;
and put $x_n=T^nx_0$, $n\ge 0$.
Denote, for simplicity reasons,
$\rho_n=\rho_n(x_0)$, $\al_n=\al_n(x_0)$, $n\ge 0$;
that is:
$\rho_n=d(x_n,x_{n+1})$, $\al_n=d(x_n,x_n)$, $n\ge 0$.
By the contractive condition 
written at ($x=x_n$, $y=x_n$), we have 
(taking (\ref{401}) into account) 
\beq \label{403}
\al_{n+1}\le \max\{\vphi(\rho_n),\al_n\},\ \ 
\forall n\ge 0.
\eeq
On the other hand, by the contractive condition 
written at ($x=x_n$, $y=x_{n+1}$), we have 
(along with (\ref{402}) above)
\beq \label{404}
\rho_{n+1}\le \max\{\vphi(\max\{\rho_n,\rho_{n+1}\}),\al_n,\al_{n+1}\},\ \ 
\forall n\ge 0.
\eeq
Now (for some $n\ge 0$), the alternative 
$\rho_{n+1}\le \vphi(\max\{\rho_n,\rho_{n+1}\})$ 
gives (as $\vphi$ is normal), $\rho_{n+1}\le \vphi(\rho_n)$;
so that, (\ref{404}) becomes
$$
\rho_{n+1}\le \max\{\vphi(\rho_n),\al_n,\al_{n+1}),\ \ 
\forall n\ge 0.
$$
This, finally combined with (\ref{403}), gives
\beq \label{405}
\rho_{n+1}\le \max\{\vphi(\rho_n),\al_n\},\ \ 
\forall n\ge 0.
\eeq
Now, by the Matthews property of the partial metric $d$,
\beq \label{406}
\al_n\le \rho_n, \ \ \forall n\ge 0.
\eeq
This, combined with the above relation, gives
$$
\rho_{n+1}\le \max\{\vphi(\rho_n),\rho_n\}\le \rho_n,\ \ 
\forall n\ge 0.
$$
The sequence $(\rho_n; n\ge 0)$ is therefore descending; 
so that,
$\rho_n\downarrow \ga$,\ \ for some $\ga\in R_+$.

{\bf (I)}
If $\ga=0$ then, by (\ref{406}), $\al_n\downarrow 0$;
so that, $(x_n; n\ge 0)$ is $e$-semi-Cauchy.
Assume now that $\ga> 0$.
By (\ref{405}),
\beq \label{407}
\ga \le \max\{\vphi(\rho_n),\al_n),\ \ \forall n\ge 0.
\eeq
As $\rho_n\downarrow \ga$ and 
$\vphi$ is right limit normal, there must be 
a sufficiently large $k=k(\ga)$ in such a way that
$\vphi(\rho_n)< \ga$, $\forall n\ge k$.
This tells us that, for $n\ge k$, 
the alternative $\ga\le \vphi(\rho_n)$ cannot hold in 
(\ref{407}); so, necessarily,
\beq \label{408}
\ga \le \al_n, \ \ \forall n\ge k.
\eeq
Combining with (\ref{406}) yields
$(\al_n; n\ge k)\downarrow \ga$; or, equivalently:
$(\al_n: n\ge 0)\searrow \ga$; hence the assertion.

{\bf (II)}
Further, we claim that $(x_n; n\ge 0)$ is $e$-Cauchy.
This, by the above obtained facts, writes
$(d(x_m,x_n); m,n\ge 0)\searrow \ga$;
or, equivalently (by (\ref{408}) above)
$(d(x_m,x_n); m,n\ge k)\downarrow \ga$.
Assume that such a property does not hold.
By Proposition \ref{p1}, there exist 
$\veps> 0$, $j(\veps)\ge k$ and a couple of 
rank-sequences 
$(m(j); j\ge k)$, $(n(j); j\ge k)$, with the
properties (\ref{206})-(\ref{209}).
Denote, for $j\ge k$,
\bit
\item[]
$t_j:=M(x_{m(j)},x_{n(j)})$, 
$s_j:=d(x_{m(j)+1},x_{n(j)+1})$,\\
$\be_j:=c(x_{m(j)},x_{n(j)})(=\max\{\al_{m(j)},\al_{n(j)}\})$.
\eit
By the quoted relations (and the $e$-semi-Cauchy property)
$$
(t_j; j\ge k)\downarrow \ga+\veps,\ 
(s_j; j\ge k)\to \ga+\veps,\ 
(\be_j; j\ge k)\downarrow \ga.
$$
On the other hand, by the contractivity condition 
at $(x=x_{m(j)}$, $y=x_{n(j)}$), 
$$
s_j\le \max\{\vphi(t_j),\be_j\},\ \ \forall j\ge k.
$$
Passing to $\limsup$ as $j\to \oo$, 
one gets (via Lemma \ref{le6} and Lemma \ref{le4})
$$
\ga+\veps\le \max\{\limsup_{t\to (\ga+\veps)+}\vphi(t),\ga\},
$$
in contradiction with the choice of $\vphi$.
Hence, the working assumption about $(x_n; n\ge 0)$ 
cannot be true; so, this sequence is $e$-Cauchy. 

{\bf (III)}
As $(X,e)$ is complete, there exists a uniquely 
determined $x^*\in X$ with $x_n\etends x^*$; 
or, equivalently (cf. Lemma \ref{le2})
$$
\lim_n d(x_n,x^*)=d(x^*,x^*)=\ga(=\lim_{m,n} d(x_m,x_n)).
$$
Denote for simplicity $\de:=d(x^*,Tx^*)$.
Clearly, $\ga\le \de$ (by the Matthews property).
Assume by contradiction that $\ga< \de$; and
let $\eta> 0$ be such that $\ga+2\eta< \de$.
As
$d(x_n,x^*)\downarrow \ga$, $d(x_n,x_{n+1})\downarrow \ga$,
there exists $k=k(\eta)$ with
\beq \label{409}
d(x_n,x^*), d(x_n,x_{n+1})< \ga+\eta< \de,\ \ 
\forall n\ge k;
\eeq
and this, by definition, yields
$M_1(x_n,x^*)=\de$, $\forall n\ge k$.
Further, by the reflexive-triangular inequality,
$$
d(x_n,Tx^*)\le d(x_n,x^*)+\de-\ga< \de+\eta;
$$
wherefrom, combining with (\ref{409}),
$$
M_2(x_n,x^*)\le (1/2)[\de+\ga+2\eta[< \de,\ \ 
\forall n\ge k.
$$
Hence, summing up,
$M(x_n,x^*)=\de$, $\forall n\ge k$.
This, along with the contractivity condition
written for $(x=x_n$, $y=x^*$), gives
$$
d(x_{n+1},Tx^*)\le \max\{\vphi(\de),\al_n,\ga\}=
\max\{\vphi(\de),\al_n\},\ \ 
\forall n\ge k.
$$
As a direct consequence, one gets 
(again by the reflexive-triangular inequality)
$$ \barr{l}
\de=d(x^*,Tx^*)\le d(x^*,x_{n+1})+d(x_{n+1},Tx^*)-\al_{n+1} \\
\le d(x^*,x_{n+1})+\max\{\vphi(\de),\al_n\}-\al_{n+1},\ \ 
\forall n\ge k.
\earr
$$
Passing to limit as $n\to \oo$ one gets 
(via Lemma \ref{le6})
$\de\le \max\{\vphi(\de),\ga\}$, contradiction.
Hence, $\ga=\de$; i.e., $x^*$ is an element
of $\Fix(T;d)$. The proof is complete.
\eproof

\section{Fixed point statement}
\setcounter{equation}{0}

Let again $X$ be a nonempty set.
Take a partial metric $d$ over it; and let $e$ stand for 
its associated metric. As before, assume that (d01) holds.

Let $T:X\to X$ be a selfmap of $X$.
Remember that, under the conditions of Theorem \ref{t1},
each $x\in X$ is a $(d;T)$-Picard point, in the sense precise there.
In particular, this tells us that
$\Fix(T;d)$ is a nonempty subset of $X$.
Note that, when $d$ is a (standard) metric, 
$\Fix(T;d)=\Fix(T)$; whence, $\Fix(T)$ is nonempty;
but, in general, this is not possible.
It is our aim in the following 
to identify -- 
in the "partial" setting --
sufficient conditions under which  
the non-emptiness of $\Fix(T)$ be still retainable.
Note that, 
by the weak sufficiency of $d$,
\beq \label{501}
z\in X,\ d(z,Tz)=d(z,z)=d(Tz,Tz) \limpl z=Tz.
\eeq
As a consequence of this, 
\beq \label{502}
z\in \Fix(T;d),\ d(z,z)=d(Tz,Tz) \limpl z\in \Fix(T).
\eeq
Hence -- 
assuming that Theorem \ref{t1} holds --
it will suffice getting points $z\in \Fix(T;d)$ 
with $d(z,z)=d(Tz,Tz)$ to conclude that 
$\Fix(T)$ is nonempty.
\sk

Technically speaking, such a conclusion is 
deductible under a stronger contractive condition 
upon $T$ than the one in Theorem \ref{t1}.
This firstly refers to the mapping $c(.,.)$ appearing there 
being substituted by $b(.,.)$; i.e., 
\bit
\item[(e01)]
$T$ is $(M;b;\vphi)$-contractive (modulo $d$).
\eit
Secondly, the function $\vphi\in \calf(R_+)$ 
in this condition must be taken so as
\bit
\item[(e02)]
$\psi:=\iota-\vphi$ is semi-coercive and $\vphi$ is limit normal.
\eit
Here, $\iota\in \calf(R_+)$ is the {\it identity} function
($\iota(t)=t$, $t\in R_+$).
It is our aim in the following to show that, under
these requirements, a positive answer to the posed
question is available.
\sk

{\bf (A)}
Note that, under such a setting, 
the working conditions of Theorem \ref{t1} hold. 
Hence, in particular, $\Fix(T;d)$ is nonempty.
For an easy reference, we shall collect some 
basic facts about this set.
These are valid even in the absence of (e01)+(e02) above.

\bprop \label{p2}
Let the selfmap $T$ of $X$ be such that 
$\Fix(T;d)\ne \es$. Then, 
\beq \label{503}
d(Tz,Tz)\le d(z,z),\ \ \forall z\in \Fix(T;d)
\eeq
\beq \label{504}
d(y,Tz)\le d(y,z),\ \ 
\forall y\in X,\ \forall z\in \Fix(T;d)
\eeq
\beq \label{505}
e(z,Tz)=d(z,z)-d(Tz,Tz),\ \ 
\forall z\in \Fix(T;d)
\eeq
\beq \label{506}
M(z,w)=d(z,w),\ \ 
\forall z,w\in \Fix(T;d).
\eeq
\eprop

\bproof
{\bf j)}
Let $z\in \Fix(T;d)$ be arbitrary fixed.
By the Matthews  property of $d$, one has
$d(Tz,Tz)\le d(z,Tz)=d(z,z)$;
and this proves (\ref{503}).

{\bf jj)}
For each $y\in X$, $z\in \Fix(T;d)$, 
we have, by the reflexive-triangular inequality,
$$
d(y,Tz)\le d(y,z)+d(z,Tz)-d(z,z)=d(y,z);
$$
hence, (\ref{504}) follows.

{\bf jjj)}
For each $z\in \Fix(T;d)$, one has, by definition,
$$
e(z,Tz)=2d(z,Tz)-d(z,z)-d(Tz,Tz)=d(z,z)-d(Tz,Tz);
$$
and, from this, (\ref{505}) is proved.

{\bf jv)}
Take some couple $z,w\in \Fix(T;d)$. We have
$$ \barr{l}
M_1(z,w)=\max\{d(z,w),d(z,Tz),d(w,Tw)\}=\\
\max\{d(z,w),d(z,z),d(w,w)\}=d(z,w);
\earr
$$
as well as, by (\ref{504}),
$$ 
M_2(z,w)=(1/2)[d(z,Tw)+d(w,Tz)]\le d(z,w).
$$
Hence, by simply combining these,  (\ref{506}) follows.
\eproof

{\bf (B)}
We may now pass to the effective part of 
our developments. 
Assume that (e01) and (e02) hold; and put
\bit
\item[(e03)]
$\theta=\inf\{d(z,z); z\in \Fix(T;d)\}$,\ 
$X(T;d)=\{z\in \Fix(T;d); d(z,z)=\theta\}$.
\eit

\bprop \label{p3}
Under the admitted conditions, we have
\beq \label{507}
\mbox{
$X(T;d)\incl \Fix(T)$;\ 
i.e.: $z\in \Fix(T;d),\ d(z,z)=\theta \limpl z=Tz$
}
\eeq
\beq \label{508}
\mbox{
$\Fix(T)$ is at most singleton; hence, so is $X(T;d)$.
}
\eeq
\eprop

\bproof
{\bf h)} 
If $\theta=0$, we are done.
For, in such a case, 
$d(z,z)=d(z,Tz)=0$ implies (by Proposition \ref{p2})
$d(Tz,Tz)=0$; whence $z=Tz$.
Assume in the following $\theta> 0$.
Taking $z$ as a starting point in Theorem \ref{t1}, one gets
(by its conclusions) that 
\beq \label{509} 
\mbox{    \btab{l}
$(\rho_n:=d(T^nz,T^{n+1}z); n\ge 0)$ is descending; 
hence $\ga:=\lim_n (\rho_n)$ exists
\etab
}
\eeq
\beq \label{510}  
\mbox{   \btab{l}
$z^*:=e-\lim_n(T^nz)$ exists, as an element of $\Fix(T;d)$, \\
with $d(z^*,z^*)=\ga$;\
hence $\theta\le d(z^*,z^*)\le d(Tz,T^2z)$.
\etab
}
\eeq
On the other hand, by the initial choice of $z$, we must have 
$\theta=d(T^nz,T^{n+1}x)$, $\forall n\ge 0$; whence $\theta=\ga$.
From the contractive condition written at ($x=z$, $y=Tz$), 
$$
\theta=d(Tz,T^2z)\le \max\{\vphi(M(z,Tz)),b(z,Tz)\}.
$$
In addition, by  (\ref{402}) and (\ref{509}),
$$
M(z,Tz)=\max\{d(z,Tz),d(Tz,T^2z)\}=d(z,Tz)=d(z,z)=\theta.
$$
Hence, by simply combining these relations,
$$
\theta\le \max\{\vphi(\theta), (1/2)[\theta+d(Tz,Tz)]\}.
$$
The alternative $\theta\le \vphi(\theta)$ is impossible.
Hence, the alternative $\theta\le (1/2)[\theta+d(Tz,Tz)]$
must be true; and then, 
$\theta\le d(Tz,Tz)\le d(z,Tz)=\theta$; wherefrom 
$d(Tz,Tz)=\theta=d(z,z)$. 
This, along with (\ref{502}), yields $z=Tz$.

{\bf hh)}
Assume that $u,v\in X$ are such that 
$u=Tu$, $v=Tv$. Then, by the Matthews property of $d$,
$$
M_1(u,v)=d(u,v),\ \ M_2(u,v)=d(u,v);
$$
hence, $M(u,v)=d(u,v)$.
By the contractive condition,
$$
d(u,v)=d(Tu,Tv)\le \max\{\vphi(d(u,v)),b(u,v)\}.
$$
If $d(u,v)\le \vphi(d(u,v))$ then
$d(u,v)=0$ (wherefrom $b(u,v)=0$); 
since this gives $e(u,v)=0$, we must have $u=v$.
If $d(u,v)\le b(u,v)$ then (combining with 
$d(u,v)\ge b(u,v)$) one has $d(u,v)=b(u,v)$.
This, by definition, yields $e(u,v)=0$;
wherefrom $u=v$.
The proof is thereby complete.
\eproof

{\bf (C)}
It therefore follows that, 
a positive answer to the posed question 
is obtainable, if we can establish that $X(T;d)$ is nonempty.
The following intermediate statement will be
useful for this purpose.

\bprop \label{p4}
Let the same conditions hold. Then,
\beq \label{511}
\forall \veps> 0, \exists z\in \Fix(T;d):\ 
d(z,z)< \theta+\veps,\  d(z,z)-d(Tz,Tz)< 2\veps.
\eeq
\eprop

\bproof
The conclusion is clear when $\theta=0$.
In fact, let $z\in \Fix(T;d)$ be such that 
$d(z,z)< \veps$;
then (by Proposition \ref{p2}),
$0\le d(Tz,Tz)\le d(z,z)< \veps< 2\veps$.
It remains now to discuss the case of $\theta> 0$.
The conclusion is again clear when $\theta\le \veps$;
because, taking $z\in \Fix(T;d)$ according to 
$d(z,z)< \theta+ \veps$, we have
$$
0\le d(Tz,Tz)\le d(z,z)< \theta+\veps\le 2\veps.
$$
So, we have to see what happens when $0< \veps< \theta$.
Assume by contradiction that the conclusion in the statement 
is not true; i.e., for some $\veps$ in $]0,\theta[$,
\bit
\item[(e04)]
$\forall z\in \Fix(T;d)$:\ $d(z,z)< \theta+\veps$ $\limpl$ 
$d(z,z)-d(Tz,Tz)\ge 2\veps$.
\eit
Note that, in such a case, one derives 
\beq \label{512}
\forall z\in \Fix(T;d):\ \ 
d(z,z)< \theta+\veps \limpl 
d(Tz,Tz)\le d(z,z)-2\veps< \theta-\veps.
\eeq
As $\vphi$ is right limit normal, there exists
some small enough $\eta$ in $]0,\veps[$, such that
\beq \label{513}
\mbox{
$\vphi(t)< \theta$,\ \ whenever $\theta\le t< \theta+\eta$.
}
\eeq
Given this $\eta> 0$, there exists some $z\in \Fix(T;d)$ with
$\theta\le d(z,z)< \theta+\eta$.
Taking $z$ as a starting point in Theorem \ref{t1}, 
the relations (\ref{509}) and (\ref{510}) are still retainable here.
From the contractive condition written at ($x=z$, $y=Tz$),
$$
\theta\le d(Tz,T^2z)\le \max\{\vphi(M(z,Tz)),b(z,Tz)\}.
$$
By (\ref{402}), the choice of $z$, and (\ref{509}), 
$$
M(z,Tz)=\max\{d(z,Tz), d(Tz,T^2z)\}=d(z,Tz)=d(z,z).
$$
Replacing in the previous relation gives
$$
\theta\le d(Tz,T^2z)\le \max\{\vphi(d(z,z)), b(z,Tz)\}.
$$
The first alternative of this, 
[$\theta\le \vphi(d(z,z))$],
is unacceptable, by the choice of $z\in \Fix(T;d)$ 
and (\ref{513}).
Moreover, the second alternative of the same,
[$\theta\le b(z,Tz)$],
yields, by the choice of $z\in \Fix(T;d)$ and (\ref{512})
$$
\theta\le (1/2)[d(z,z)+d(Tz,Tz)]< 
(1/2)[\theta+\eta+\theta-\eta]=\theta;
$$
again a contradiction.
So, (e04) cannot be true; and conclusion follows.
\eproof

Having these precise, 
we may now answer the  posed question.

\bprop \label{p5}
Let (e01) and (e02)
be in use. Then, $X(T;d)$ is nonempty.
\eprop

\bproof
Let $(\mu_n; n\ge 0)$ be a strictly descending sequence 
in $R_+^0$, with $\mu_n\to 0$.
(For example, one may take $\mu_n=2^{-n}$, $n\ge 0$).
By Proposition \ref{p4}, there exists, for each $n\ge 0$, 
some $z_n\in \Fix(T;d)$ with
\beq \label{514}
d(z_n,z_n)< \theta+\mu_n,\ \ d(z_n,z_n)-d(Tz_n,Tz_n)< 2\mu_n.
\eeq
Let $m,n\ge 0$ be arbitrary for the moment. 
We have (by the reflexive-triangular property 
and $(z_n; n\ge 0)\incl \Fix(T;d)$)
$$  \barr{l}
d(z_m,z_n)\le d(z_m,Tz_m)+d(Tz_m,z_n)-d(Tz_m,Tz_m)= \\
d(z_m,z_m)-d(Tz_m,Tz_m)+d(Tz_m,z_n)< 2\mu_m+d(Tz_m,z_n).
\earr
$$
The distance in the right hand side may be evaluated as
$$ \barr{l}
d(Tz_m,z_n)\le d(Tz_m,Tz_n)+d(Tz_n,z_n)-d(Tz_n,Tz_n)= \\
d(Tz_m,Tz_n)+d(z_n,z_n)-d(Tz_n,Tz_n)< 2\mu_n+d(Tz_m,Tz_n).
\earr
$$
So, replacing in the previous one, we get 
$$
d(z_m,z_n)\le 2[\mu_m+\mu_n]+d(Tz_m,Tz_n),\ \ \forall m,n\ge 0;
$$
and this, combined  with the contractive condition, yields
(via Proposition \ref{p2})
\beq \label{515}
d(z_m,z_n)\le 2(\mu_m+\mu_n)+
\max\{\vphi(d(z_m,z_n)),b(z_m,z_n)\},\ \  
\forall m,n\ge 0.
\eeq

{\bf Part 1}.
Denote, for simplicity
\bit
\item[]
$\al_n:=d(z_n,z_n)$,\ $\rho_n:=d(z_n,z_{n+1})$,\
$\be_n:=b(z_n,z_{n+1})$,\ \ $n\ge 0$.
\eit
By (\ref{514}), $(\al_n; n\ge 0) \downarrow \theta$;
and this, by definition, yields 
$(\beta_n; n\ge 0) \downarrow \theta$.
On the other hand, by the Matthews property of $d$, 
$$ \mbox{
$\theta\le \al_n\le \rho_n$, $\forall n\ge 0$;\ 
hence $\theta\le \sig:=\limsup_n(\rho_n)$.
}
$$
We claim that $\sig=\theta$; 
or, equivalently (see above) $\rho_n \downarrow \theta$.
Suppose by contradiction that $\sig> \theta$;
note that $\sig\le \oo$.
As a direct consequence of (\ref{515}), 
\beq \label{516}
\rho_n\le 2[\mu_n+\mu_{n+1}]+ 
\max\{\vphi(\rho_n),\be_n\},\ \ \forall n\ge 0.
\eeq
On the other hand, by the definition of $\sig$, 
there must be a rank-sequence $(n(j); j\ge 0)$ with
$n(j)\to \oo$ and $\rho_{n(j)}\to \sig$.
Let $\veps> 0$ be such that 
$\theta+ \veps< \theta+ 2\veps< \sig$.
As $\be_n\downarrow \theta$, $\mu_n\downarrow 0$,
there must be some $j(\veps)$, in such a way that
$$
\rho_{n(j)}> \theta+2\veps,\
\be_{n(j)}< \theta+ \veps,\
\mu_{n(j)}< \veps/4,\ \ \forall j\ge j(\veps);
$$
whence, it is the case that
$$
\rho_{n(j)}> 2[\mu_{n(j)}+\mu_{n(j)+1}]+\be_{n(j)},\  
\forall j\ge j(\sig).
$$
This, along with (\ref{516}), shows us that 
we must have (as alternative)
\beq \label{517}
\rho_{n(j)}\le 2[\mu_{n(j)}+\mu_{n(j)+1}]+\vphi(\rho_{n(j)}),\ \ 
\forall j\ge j(\veps).
\eeq

{\bf Sub-case 1a}. 
Suppose that $\sig=\oo$.
By the introduced in (e02) notation, 
$$
\psi(\rho_{n(j)})\le 2[\mu_{n(j)}+\mu_{n(j)+1}],\ \ 
\forall n\ge 0.
$$
Then, passing to limit as $j\to \oo$, one gets
$\lim_j \psi(\rho_{(n(j)})=0$; in contradiction 
with the property imposed upon $\psi$ and Lemma \ref{le5}.

{\bf Sub-case 1b}. 
Suppose that $\sig< \oo$.
Passing to limit as $j\to \oo$ in (\ref{517}),
one gets
$\sig\le \limsup_{t\to \sig} \vphi(t)$; 
again a contradiction.

Summing up, $(\rho_n; n\ge 0)\downarrow \theta$;
wherefrom, $(z_n; n\ge 0)$ is $e$-semi-Cauchy.

{\bf Part 2}.
We show that $(z_n; n\ge 0)$ is $e$-Cauchy.
Suppose that this is not true.
By Proposition \ref{p1},
there exist then $\veps> 0$, $j(\veps)\ge 0$, and a couple
of rank-sequences $(m(j); j\ge 0)$, $(n(j); j\ge 0)$ with
the properties (\ref{206})-(\ref{209}).
Further, by (\ref{515}), 
\beq \label{518}
\barr{l}
d(z_{m(j)},z_{n(j)})\le 2[\mu_{m(j)}+\mu_{n(j)}]+ \\
\max\{\vphi(d(z_{m(j)},z_{n(j)})),b(z_{m(j)},z_{n(j)})\},\ \  
\forall j\ge 0.
\earr
\eeq
Passing to limit as $j\to \oo$, one gets, via 
(\ref{208}), Lemma \ref{le6} and Lemma \ref{le4},
$$
\theta+\veps\le \max\{\limsup_{t\to \theta+\veps}\vphi(t),\theta\};
$$
impossible by the choice of $\vphi$.
Hence the conclusion.

{\bf Part 3}.
As $(X,e)$ is complete, there exists 
$z\in X$ with $z_n\etends z$; i.e.,
$$
d(z,z)=\lim_nd(z_n,z)=\lim_{m,n} d(z_m,z_n)=\theta.
$$
We claim that $z\in \Fix(T;d)$; this, by the relation above,
amounts to $d(z,Tz)=\theta$.
For the moment, we have (cf. Proposition \ref{p2})
$$
d(Tz_n,Tz_n)\le d(z_n,z_n)\le d(z,z_n),\ \  \forall n\ge 0;
$$
hence
$$
\de_n:=d(z,z_n)-d(Tz_n,Tz_n)\ge 0,\ \ \forall n\ge 0.
$$
In addition, by the choice of $(z_n; n\ge 0)$,
$$
\de_n=d(z,z_n)-d(z_n,z_n)+d(z_n,z_n)-d(Tz_n,Tz_n)\le 
d(z,z_n)-\theta+2\mu_n,\ \ \forall n\ge 0;
$$
and this gives $\de_n\to 0$ as $n\to \oo$.
On the other hand, by the reflexive-triangular inequality
and Proposition \ref{p2}, one has, for all $n$,
$$
d(z,Tz)\le d(z,z_n)+d(z_n,Tz)-d(z_n,z_n),
$$
$$ \barr{l}
d(z_n,Tz)\le d(z_n,Tz_n)+d(Tz_n,Tz)-d(Tz_n,Tz_n)\le \\
d(z_n,z_n)+d(Tz_n,Tz)-d(Tz_n,Tz_n);
\earr
$$
so, by combining these,
\beq \label{519}
d(z,Tz)\le \de_n+d(Tz_n,Tz),\ \ 
\forall n\ge 0.
\eeq
Now, $d(z,Tz)\ge d(z,z)=\theta$.
Assume by contradiction that $d(z,Tz)> \theta$;
and let $\veps> 0$ be such that
$\theta+ 2\veps< d(z,Tz)$.
There exists $k=k(\veps)$ such that, for all $n\ge k$,
$$
d(z_n,z)< \theta+\veps< d(z,Tz),
$$
By the contractive condition,
$$
d(Tz_n,Tz)\le \max\{\vphi(M(z_n,z)),b(z_n,z)\},\ \ 
\forall n\ge 0.
$$
To evaluate the right member, we have 
(from Proposition \ref{p2} and the preceding relation), 
for all $n\ge k$,
$$ \barr{l}
M_1(z_n,z)=\max\{d(z_n,z),d(z_n,Tz_n),d(z,Tz)\}=\\
\max\{d(z_n,z),d(z_n,z_n),d(z,Tz)\}=
\max\{d(z_n,z),d(z,Tz)\}=d(z,Tz),
\earr
$$
$$ \barr{l}
M_2(z_n,z)=(1/2)[d(z_n,Tz)+d(Tz_n,z)]\le \\
(1/2)[d(z_n,z)+d(z,Tz)-d(z,z)+d(z_n,z)]=\\
(1/2)[2d(z_n,z)-d(z,z)+d(z,Tz)]< \\
(1/2)[\theta+2\veps+d(z,Tz)]< d(z,Tz).
\earr
$$
This yields $M(z_n,z)=d(z,Tz)$, $\forall n\ge k$; 
so that, the contractive condition becomes
$$
d(Tz_n,Tz)\le \max\{\vphi(d(z,Tz)),b(z_n,z)\},\ \ 
\forall n\ge k;
$$
wherefrom, replacing in (\ref{519}),
\beq \label{520}
d(z,Tz)\le \de_n+\max\{\vphi(d(z,Tz)),b(z_n,z)\},\ \ 
\forall n\ge k.
\eeq
Note that
$$
\mbox{
$b(z_n,z)=(1/2)[d(z_n,z_n)+d(z,z)]\to \theta$,\ \ 
as $n\to \oo$.
}
$$
So, passing to limit as $n\to \oo$ in (\ref{520}), 
gives, via Lemma \ref{le6}, 
$$
d(z,Tz)\le \max\{\vphi(d(z,Tz)),\theta\};
$$
contradiction. 
Hence, the working assumption about $d(z,Tz)$ 
cannot be accepted; wherefrom $d(z,Tz)=\theta$.
This completes the argument.
\eproof

{\bf (D)}
Now, as a direct consequence of all these, 
we have the second main result of this exposition.

\btheorem \label{t2}
Suppose that (e01) and (e02) hold. Then,
$\Fix(T)$ is a nonempty singleton;
i.e.: $T$ has a unique fixed point in $X$.
\etheorem

In particular, assume $\vphi$ is linear;
i.e., $\vphi(t)=\al t$, $t\in R_+$, for some $\al$ in $[0,1[$.
Then, Theorem \ref{t2} is just the main result in 
Ili\'{c} et al \cite{ilic-pavlovic-rakocevic-2012}.
But, we must say that the way of proving it
(via Theorem \ref{t1}) is different from the one proposed in
that paper.
Further aspects may be found in 
Chi et al \cite{chi-karapinar-thanh-2012}.


\end{document}